\documentclass[12pt,draftcls,onecolumn]{IEEEtran}
\usepackage{fancyhdr}

\usepackage{amsfonts}
\usepackage{mathrsfs}
\usepackage{}
\usepackage{colortbl}
\usepackage{amssymb,amsmath,graphicx}

\usepackage{caption}
\usepackage{subfigure}
\usepackage{float}

\newtheorem{Remark}{\bf Remark}[section]

\newtheorem{Corollary}{\bf Corollary}[section]

\newtheorem{Definition}{\bf Definition}[section]

\newenvironment{Proof}{\noindent{\em Proof:\/}}{\hfill $\Box$\par}

\newtheorem{Theorem}{\bf Theorem}[section]

\newtheorem{Lemma}{\bf Lemma}[section]

\newtheorem{Property}{Property}[section]
\newtheorem{Assumption}{\bf Assumption}[section]

\newcommand{\R}{{\mathbb R}}

\newcommand{\EQQ}{\begin{eqnarray*}}
\newcommand{\ENN}{\end{eqnarray*}}
\newcommand{\EQ}{\begin{eqnarray}}
\newcommand{\EN}{\end{eqnarray}}
\newcommand{\bd}{\begin{Definition}
\newcommand{\mbox{col}}{\mbox{col}}
\begin{rm} }
\newcommand{\ed}{ \end{rm}
\end{Definition} }

\begin{document}

\title{\LARGE \bf
Cooperative Output Regulation of Discrete-Time Linear Time-Delay Multi-agent Systems
}



\author{Yamin~Yan~and~Jie~Huang
\thanks{This work has been supported by the Research Grants Council of the Hong Kong Special Administration
Region under grant No. 412813.}
\thanks{Yamin Yan and Jie Huang  are with Department of Mechanical and Automation
Engineering, The Chinese University of Hong Kong, Shatin, N.T., Hong
Kong. E-mail: ymyan@mae.cuhk.edu.hk, jhuang@mae.cuhk.edu.hk}
}

\maketitle

\thispagestyle{fancy}
\fancyhead{}
\lhead{}
\lfoot{}
\cfoot{}
\rfoot{}
\renewcommand{\headrulewidth}{0pt}
\renewcommand{\footrulewidth}{0pt}

\begin{abstract}
In this paper, we study the cooperative output regulation problem
for the discrete linear time-delay multi-agent systems by distributed observer approach. In contrast with the same problem for continuous-time  linear time-delay multi-agent systems, the problem has two new features. First, in the presence of
 time-delay, the  regulator equations for discrete-time linear systems are different from those for continuous-time linear systems.
Second,  under the standard assumption on the connectivity of the communication graph,
a distributed observer for  any continuous-time linear leader system always exists. However, this is not the case for the discrete counterpart of the distributed observer for the continuing-time systems. We will propose another type of discrete distributed observer. It turns out that such an observer always exists under some mild assumptions.  Using this result, we further present
 the solvability conditions of the problem by distributed dynamic output feedback control law.
\end{abstract}

\section{Introduction}

The cooperative output regulation problem aims to design a distributed control law for a multi-agent system to drive the tracking error of each follower to the origin asymptotically while rejecting a class of external disturbances. It can be viewed as an extension of the classical output regulation problem from a single system to a multi-agent system. It can also be viewed as an extension of the leader-following consensus problem \cite{jadbabaie2003, olfatisaber2004, Hu1} in that it not only handles the asymptotically tracking problem, but also the disturbance
rejection problem.

The classical output regulation problem for a single system was thoroughly studied for both continuous-time linear systems \cite{Da,Fr,FW} and discrete-time linear systems \cite{DW, J.H}.
The key for the solvability of the problem is to obtain the solution of a set of matrix equations called the  regulator equations. Since the  regulator equations for both continuous-time linear systems and  discrete-time linear systems are the same, the extension of the results  on the output regulation problem from  continuous-time linear systems to discrete-time linear systems is quite straightforward.

The cooperative output regulation problem was first studied for continuous-time linear multi-agent systems  in \cite{Su1,Su5}. What makes this problem interesting is that the control law has to satisfy certain communication constraints described by a digraph.  Such a control law is called a distributed control law.
The core of the approach in \cite{Su1,Su5} is the employment of a so-called distributed observer which is a dynamic compensator capable of providing the estimation of the leader's signal to each follower without violating the communication constraints.
On the basis of the distributed observer, a distributed controller satisfying all the communication constraints can be synthesized to solve the cooperative output regulation problem.

Recently, the cooperative output regulation problem for continuous-time linear multi-agent systems with time-delay was also studied in \cite{CCC14}.
It turns out that the distributed observer approach also played a key role in constructing a distributed control law for solving the problem for time-delay systems.

In this paper, we will consider the cooperative output regulation problem for discrete-time linear multi-agent systems with time-delay by a distributed control law.
In contrast with the same problem for continuous-time  linear systems, the problem in this paper has two new features. First, as pointed out in \cite{CCC15},
for linear systems with time-delay, the  regulator equations for discrete-time linear systems are different from those for continuous-time linear systems.
Second and more intriguingly, under the standard assumption on the connectivity of the communication graph,
a distributed observer for any continuous-time linear leader system always exists by having the observer gain sufficiently large because a sufficient large observer gain can place the eigenvalues
far left  of the complex plane.
However, for discrete-time systems, since the stability region is the unit circle, a large observer gain may destabilize a system. In fact, it will be pointed out in Section \ref{dob} that, the discrete counterpart of the distributed observer
in \cite{Su1} may never exist when the eigenvalues of the system matrix of the leader system have different sign. Thus, instead of using the discrete counterpart of the distributed observer
in \cite{Su1}, we propose another discrete distributed observer and give a thorough study on the stability of this observer in Section \ref{dob}. It will be shown that the existence condition
of this new type of observer is much  milder than the discrete counterpart of the distributed observer in \cite{Su1}. In particular, this observer always exists for a marginally stable leader systems provided that the communication graph is connected.

The consensus problem of discrete-time systems  is studied in \cite{you2011, you2013} for delay-free case,  and in \cite{hulin, xiecheng} for time-delay case. However, references \cite{hulin, xiecheng}
focus on systems in chain-integrator form. In contrast, this paper deals with a general linear system and addresses both asymptotic tracking and disturbance rejection.

The rest of this paper is organized as follows. Section \ref{cpfp}
formulates the problem and lists some assumptions. Section \ref{dob} studies  the stability property of the discrete distributed observer. Section \ref{cmr}
presents our main results. An example is used to illustrate our design in Section \ref{ste}. Finally we close the paper with  some concluding remarks in Section \ref{cc}.

\textbf{Notation.} $\otimes$ denotes the Kronecker product of matrices. $\sigma(A)$ denotes the spectrum of a square matrix $A$. $\mathbb{Z^+}$ denotes the set of nonnegative integers. $\mathbf{1}_N$ denotes an $N\times1$ column vector whose elements are all $1$.
For $X_i \in \R^{n_i}$, $i = 1, \dots, m$,
$\mbox{col} (X_1, \dots, X_m) = [X^T_1, \dots, X^T_m]^T$.
For some nonnegative integer $r$, $I[-r,0]$ denotes the set of integers $\{-r,-r+1,\cdots,0\}$ , and $\mathcal{C}(I[-r,0],\R^{n})$ denotes the set of functions mapping the integer set $I[-r,0]$ into $\R^{n}$.

\section{Problem Formulation and
Preliminaries}\label{cpfp}

Consider the discrete time-delay multi-agent systems as follows:
\begin{equation}\label{sys1}
\begin{split}
x_i(t+1)=&A_ix_i(t)+\sum_{l=0}^hB_{il}u_i(t-r_l)+E_{xi}x_0(t)+E_{wi}w_i(t),\\
y_{mi}(t)=&\bar{C}_ix_i(t)+\sum_{l=0}^h\bar{D}_{il}u_i(t-r_l)+\bar{F}_{xi}x_0(t)+\bar{F}_{wi}w_i(t),\\
e_i(t)=&C_ix_i(t)+\sum_{l=0}^hD_iu_i(t-r_l)+F_{xi}x_0(t)+F_{wi}w_i(t),\\
 & t\in \mathbb{Z}^+, ~i=1,\cdots, N\\
\end{split}
\end{equation}
where $x_i\in \R^{n_i}$ is the state, $u_i\in \R^{m_i}$ the input, $e_i\in \R^{p_i}$ the error output, $y_{mi}\in \R^{p_{mi}}$ the measurement output, $x_0\in \R^{q}$ the measurable exogenous signal such as the reference input to be tracked and $w_i\in \R^{s_i}$ the unmeasurable exogenous signal  such as external disturbances to each subsystem. $r_l\in \mathbb{Z^+}$, $l=0,\cdots,h$, satisfying $0=r_0<r_1<\cdots<r_h=r<\infty$.  We assume that $x_0$ is generated by the exosystem of the following form
\begin{equation}\label{exo1}
\begin{split}
x_0(t+1)&=S_0x_0(t), ~ t\in \mathbb{Z}^+,\\
\end{split}
\end{equation}
where $S_0\in \R^{q\times q}$ is a constant matrix.

Also, it is assumed that, for $i=1,\cdots,N$, $w_i$ is generated by a linear system as follows
\begin{equation}\label{distur}
\begin{split}
w_i(t+1)&=Q_iw_i(t), t \in \mathbb{Z}^+,\\
\end{split}
\end{equation}
with $Q_i\in\R^{s_i\times s_i}.$

Like in \cite{CCC14}, we can view the system composed of \eqref{sys1} and \eqref{exo1} a multi-agent systems with the exosystem \eqref{exo1} being the leader and the $N$ subsystems of \eqref{sys1} being the $N$ followers. When $N=1$, this multi-agent problem reduces to a single system as studied in \cite{CCC15}. Associated with this multi-agent system, we can define a digraph
$\mathcal{\bar{G}}=(\mathcal{\bar{V}},\mathcal{\bar{E}})$\footnote{See Appendix A  for
a summary of graph.}
where $\mathcal{\bar{V}}=\{0,1,\dots,N\}$ with the node 0 associated with
the exosystem \eqref{exo1} and all the other nodes
associated with the $N$ subsystems of \eqref{sys1}, and
$(i,j) \in \mathcal{\bar{E}}$ if and only if the control $u_j$ can access the measurable output $y_{mi}$  of the subsystem $i$.
We denote the adjacent matrix of  $\mathcal{\bar{G}}$ by $\bar{\mathcal{A}}=[a_{ij}]\in R^{(N+1)\times (N+1)}, i,j=0,\cdots,N$,  and
the neighbor set of the node $i$ by $\mathcal{\bar{N}}_i $.

We will consider the output feedback control law of the following form:
\begin{equation}\label{ctr22}
\begin{split}
u_i(t)=& k_i (\zeta_i(t))\\
\zeta_i(t+1)=& g_i \big(\zeta_i(t), \zeta_i(t-r_1), \cdots, \zeta_i(t-r_l),  y_{mi}(t), \\& \zeta_j(t), \zeta_j(t-r_1), \cdots, \zeta_j(t-r_l), y_{mj}(t), j\in\bar{\mathcal{N}}_i\big)\\
i&=1,\cdots,N
\end{split}
\end{equation}
where $y_{m0} = x_0$,  $\zeta_{i}\in \R^{l_i}$ for some integer $l_i$, and $k_i$ and $g_i$ are linear functions of their arguments.
It can be seen that, for each $i = 1, \cdots, N$, $j = 0, 1, \cdots, N$, $u_i$ of (\ref{ctr22}) depends on $y_{mj}$ only if the agent $j$ is a neighbor of the agent $i$. Thus,
the control law (\ref{ctr22}) is a distributed control law. The specific  control law will be given in Section \ref{cmr}.
%
%

Now we describe the cooperative output regulation problem as
follows.

\begin{Definition}\label{de2.1}
Given the multi-agent system composed of \eqref{sys1}, \eqref{exo1}, \eqref{distur} and the digraph $\bar{\mathcal{G}}$, find a
distributed control law of the form \eqref{ctr22} such that the following properties hold.
\begin{Property}\label{per21}
The origin of the closed-loop system is exponentially stable when $x_0=0$ and $w_i = 0$, $i = 1, \cdots, N$.
\end{Property}
\begin{Property}\label{per22}
For any initial condition $x_{00}  \in \mathcal{C}(I[-r,0],\R^{q})$, and,  $x_{i0}  \in \mathcal{C}(I[-r,0],\R^{n_i})$, $\zeta_{i0}  \in \mathcal{C}(I[-r,0],\R^{l_i})$, and $w_{i0} \in \R^{s_i}$,
$i = 1, \cdots, N$,  the trajectory of the closed-loop system satisfy
\end{Property}
\begin{equation}\label{}
\lim_{t \rightarrow \infty}e_i(t)=0,\ i=1,\dots,N.
\end{equation}
\end{Definition}

In order to solve the problem, we need some standard
assumptions  as follows.


\begin{Assumption}\label{Ass2.2}
There exist matrices $K_{1i} \in \R^{m_i \times n_i},\ i=1,\dots,N$, such that the system $x_i(t+1) =A_{i}x_i(t) + \sum_{l=0}^{h}B_{il}  K_{1i} x_i(t-r_l)$ is exponentially stable.
\end{Assumption}

\begin{Assumption}\label{Ass2.3}
 $ \bigg(\left[\begin{array}{cc}
                                            A_i&E_{wi}\\
                                            0_{s_i\times n_i}&Q_i\\
                                            \end{array}\right], \left[\begin{array}{cc}
                                                                \bar{C}_i&\bar{F}_{wi}
                                                                \end{array}\right]\bigg)$, $i=1,\dots,N$, are detectable.
\end{Assumption}
\begin{Assumption}\label{Ass2.4}
The matrix equations
\begin{equation}\label{re23}
\begin{split}
X_i \bar{S}_i^{r+1} &=  A_{i} X_i \bar{S}_i^r + \sum_{l=0}^{h} B_{il} U_i \bar{S}_i^{r-r_l} + E_i\bar{S}^r_i \\
0     &=  C_{i} X_i\bar{S}_i^r  + \sum_{l=0}^{h} D_{il} U_i  \bar{S}_i^{r-r_l} + F_i\bar{S}^r_i,\ i=1,\dots,N \\
\end{split}
\end{equation}
where $\bar{S_i}=\bigg[\begin{array}{cc}
                            S_0&0_{q\times s_i}\\
                            0_{s_i \times q}&Q_i\\
                            \end{array}\bigg]$, $E_i=[\begin{array}{cc}
                                                        E_{xi}&E_{wi}\\
                                                        \end{array}]$ and $F_i=[\begin{array}{cc}
                                                                        F_{xi}&F_{wi}\\
                                                                        \end{array}]$, have solution pairs $(X_i,U_i)$, respectively.
\end{Assumption}

\begin{Assumption}\label{Ass2.5}
The  graph $\bar{\mathcal{G}}$ is static and connected, i.e. $\bar{\mathcal{E}}$ is time-invariant and every node $i = 1,\cdots,N$ is reachable from the node 0 in $\bar{\mathcal{G}}$.
\end{Assumption}

\begin{Assumption}\label{Ass2.7}
All the eigenvalues of $S_0$ have modulus equal to or smaller than $1$.
\end{Assumption}


\begin{Remark}\label{rem23.1}
Assumption \ref{Ass2.2} can be viewed as a stabilizability condition for time-delay systems. By Remark 2.1 of \cite{CCC15},  Assumption \ref{Ass2.2} is satisfied if and only if
there exists a matrix $K_{1i} \in \R^{m_i\times n_i}$ such that all the roots of the following polynomial
\begin{equation}\label{eqcp}
\begin{split}
\Delta _i(z)= det{ \left(z I_{n_i} - A_i- \sum_{l=0}^{h}  B_{il} K_{1i} z^{- r_l} \right) }
\end{split}
\end{equation}
have modulus smaller than $1$.
As for Assumption \ref{Ass2.4}, Equations \eqref{re23} are called the discrete regulator equations  which are different from those of continuous-time systems \cite{CCC14}. From \cite{CCC15}, Assumption \ref{Ass2.4} is satisfied  if
for all $z \in \sigma(\bar{S}_i)$,
\begin{equation}\label{tzero2}
rank \left(
                                      \begin{array}{cc}
                                         z^rA_i -z^{r+1} I_{n_i} & \sum_{l=0}^{h} B_{il} z^{r-r_l} \\
                                         z^rC_i  & \sum_{l=0}^{h} D_{il} z^{r-r_l}\\
                                      \end{array}
                                    \right)=n_i+p_i.
\end{equation}
Assumption \ref{Ass2.5} is a standard assumption in the literature of the cooperative control of multi-agent systems under static networks.
From the adjacent matrix  $\mathcal{\bar{A}}$ of the digraph $\mathcal{\bar{G}}$, we can define a square matrix
$H=[h_{ij}]\in \R^{N\times N}$ with $h_{ii}=\sum_{j=0}^Na_{ij}$ and $h_{ij}=-a_{ij}$, for any $i\neq j$. Under Assumption \ref{Ass2.5}, $-H$ is Hurwitz \cite{Hu1}. Thus, if the digraph $\mathcal{\bar{G}}$ is also undirected, $H$ is symmetric and positive definite.

Assumption \ref{Ass2.7}  is not restrictive since it is satisfied by a large class of signals such as the step function, ramp function, and sinusoidal function.
\end{Remark}

\vspace{0.3cm}

\section{Discrete Distributed Observer}\label{dob}
Recall from \cite{Su1} that, given a continuous-time linear leader system of the form $\dot{x}_0 = S_0 x_0$, we can define a dynamic compensator of the following form

\begin{equation}\label{cob}
\begin{split}
\dot{\eta}_i(t)&=S_0\eta_i(t)+\mu\bigg(\sum_{j\in\bar{\mathcal{N}}_i}a_{ij}\big(\eta_j(t)-\eta_i(t)\big)\bigg),~ i =1,\cdots,N\\
\end{split}
\end{equation}
where $\eta_0(t)=x_0(t)$, and $\mu$ is a real number called observer gain. It can be seen that $\dot{\eta}_i(t)$ depends on $\eta_j(t)$ for $j=0,\cdots,N$ and $j\neq i$ iff $a_{ij}\neq 0$.
Let $\eta(t)= \mbox{col} (\eta_1(t),\cdots, \eta_N(t))$, $\hat{x}_0(t)=\mathbf{1}_N\otimes x_0(t)$, and $\tilde{\eta}(t)=\eta(t)-\hat{x}_0(t)$. Then it can be verified that  $\tilde{\eta}(t)$ satisfies
\begin{equation}\label{obtd}
\begin{split}
\dot{\tilde{\eta}} (t)&=\big( (I_N\otimes S_0)-\mu(H\otimes I_q)\big)\tilde{\eta}(t)\\
\end{split}
\end{equation}
Thus, if there exists some $\mu$ such that the matrix $\big( (I_N\otimes S_0)-\mu(H\otimes I_q)$ is Hurwitz, then, for any $x_0
(0)$, and $\eta_i (0)$, $ i = 1, \dots, N$, we have \EQ \label{p0}
\lim_{t \rightarrow \infty} (\eta_i(t) - x_0 (t)) = 0. \EN
Thus, we call the system \eqref{cob} a distributed observer of the leader system $\dot{x}_0 = S_0 x_0$ if and only if the system \eqref{obtd} is asymptotically stable.
It was shown in  \cite{Su1} that, under Assumption \ref{Ass2.5}, for any matrix $S_0$, there exists sufficiently large $\mu$ such that the matrix $\big( (I_N\otimes S_0)-\mu(H\otimes I_q)$ is Hurwitz.

The discrete counterpart of (\ref{cob}) is as follows:
\begin{equation}\label{ob}
\begin{split}
\eta_i(t+1)&=S_0\eta_i(t)+\mu\bigg(\sum_{j\in\bar{\mathcal{N}}_i}a_{ij}\big(\eta_j(t)-\eta_i(t)\big)\bigg),~i =1,\cdots,N\\
\end{split}
\end{equation}
where $\eta_0(t)=x_0(t)$, and $\mu$ is a real number. With $\eta(t)= \mbox{col} (\eta_1(t),\cdots, \eta_N(t))$, $\hat{x}_0(t)=\mathbf{1}_N\otimes x_0(t)$, and $\tilde{\eta}(t)=\eta(t)-\hat{x}_0(t)$, (\ref{ob}) can be put in the following form:
\begin{equation}\label{obt}
\begin{split}
\tilde{\eta} (t+1)&=\big( (I_N\otimes S_0)-\mu(H\otimes I_q)\big)\tilde{\eta}(t)\\
\end{split}
\end{equation}

Like \eqref{cob}, (\ref{ob}) is a discrete distributed observer of the leader system  if and only if there exists some $\mu$ such that the matrix $\big( (I_N\otimes S_0)-\mu(H\otimes I_q)$ is Schur.
Nevertheless, the stability property of the matrix $\big( (I_N\otimes S_0)-\mu(H\otimes I_q)$ is much more complicated than the continuous time case. For example consider a simple case
where $S_0 = \mbox{diag} [1, -1]$, $H$ is any symmetric positive matrix with eigenvalues $0 < \lambda_1 \leq \cdots \leq  \lambda_N$. Then the eigenvalues of
the matrix $\big( (I_N\otimes S_0)-\mu(H\otimes I_q)$ are given by $\{ \pm 1 - \mu \lambda_l, \; l = 1, \cdots, N \}$. Thus,   the matrix $\big( (I_N\otimes S_0)-\mu(H\otimes I_q)$ is Schur if and only if  $|\pm 1 - \mu \lambda_l| < 1$, $l= 1, \cdots, N$, and only if $\mu > 0$ and $\mu <0$. As a result, there exists no $\mu$ to make the matrix $\big( (I_N\otimes S_0)-\mu(H\otimes I_q)$  Schur.

Thus, instead of (\ref{ob}), we propose the following candidate for the discrete distributed observer of the leader system \eqref{exo1}:
\begin{equation}\label{ob2}
\begin{split}
\eta_i(t+1)&=S_0\eta_i(t)+\mu S_0 \bigg(\sum_{j\in\bar{\mathcal{N}}_i}a_{ij}\big(\eta_j(t)-\eta_i(t)\big)\bigg),~ i =1,\cdots,N\\
\end{split}
\end{equation}
whose compact form is as follows:

\begin{equation}\label{obt2}
\begin{split}
\tilde{\eta} (t+1)&=\big( (I_N\otimes S_0)-\mu(H\otimes S_0)\big)\tilde{\eta}(t)\\
\end{split}
\end{equation}

To give a detailed study on the stability for the system of the form (\ref{obt}),
denote the eigenvalues of $S_0$ by $\{\lambda_1, \cdots, \lambda_q\}$ where $0\leq |\lambda_1| \leq \cdots \leq |\lambda_q|$, and the eigenvalues of $H$ by $\{a_l\pm jb_l\}$, where $b_l=0$ when $1\leq l \leq N_1$ with $0 \leq N_1 \leq N$ and $b_l\neq 0$ when $(N_1+1)\leq l\leq N_2$ where $N_1+2(N_2-N_1)=N$.
%

Then, we have the following result.

\begin{Lemma}\label{Theo2}
Under Assumption \ref{Ass2.5},  the matrix $\big( (I_N\otimes S_0)-\mu(H\otimes S_0)\big)$ is Schur for some real $\mu$ iff
\EQ \label{con1}
{a_l^2+b_l^2} >  {b^2_{l}}  |\lambda_q|^2,~ l = 1, \cdots, N
\EN
and,
\EQ \label{mu1}
\max_{l = 1, \cdots, N} \bigg \{ \frac{ a_l -\sqrt{\Delta_l }}{a_l^2+b_l^2}\bigg  \} < \min_{l = 1, \cdots, N} \bigg \{ \frac{ a_l  + \sqrt{\Delta_l}}{a_l^2+b_l^2} \bigg \}.
\EN
where $\Delta_l = \frac{(a_l^2+b_l^2)}{ |\lambda_q|^2 } -  b^2_{l}$.

If the conditions (\ref{con1}) and (\ref{mu1}) are satisfied, the matrix $\big( (I_N\otimes S_0)-\mu(H\otimes S_0)\big)$ is Schur for all $\mu$ satisfying
\EQ \label{conmu}
 \max_{l = 1, \cdots, N} \bigg \{ \frac{ a_l -\sqrt{\Delta_l }}{a_l^2+b_l^2} \bigg \} < \mu  < \min_{l = 1, \cdots, N} \bigg \{ \frac{ a_l  + \sqrt{\Delta_l}}{a_l^2+b_l^2} \bigg \}.
\EN

%

\begin{Proof}
By the property of Kronneker product, it can be verified that the eigenvalues of $(I_N\otimes S_0-\mu(H\otimes S_0))$ are
\EQ \label{eig}
\bigg \{ (1-\mu (a_l\pm j b_l))\lambda_k, ~ k=1,\cdots,q, l=1,\cdots,N \bigg \}
\EN
Thus the matrix  $(I_N\otimes S_0-\mu(H\otimes S_0))$ is Schur if and only if, for  $l=1,\cdots,N$,
\EQQ
|1-\mu (a_l\pm j b_l)|<\frac{1}{|\lambda_q|}
\ENN
and if and only if, for  $l=1,\cdots,N$,
\EQ \label{eqleq1}
(1-\mu a_l)^2+(\mu b_l)^2<\frac{1}{|\lambda_q|^2}.
\EN
or
\EQ \label{minus}
(a_l^2+b_l^2)\mu^2-2 a_l \mu + 1 - \frac{1}{|\lambda_q|^2}<0
\EN
For each  $l=1,\cdots,N$, let $p_{l}  (\mu) = (a_l^2+b_l^2)\mu^2-2 a_l \mu+ 1 - \frac{1}{|\lambda_q|^2}$.
Then it can be verified that $p_{l} (\mu) $ has two distinct real roots if and only if $\Delta_l = \frac{(a_l^2+b_l^2)}{ |\lambda_q|^2 } -  b^2_{l}>0$.
Under Assumption 2.4, $a_l^2+b_l^2>0$. Thus,  for all $\mu \in \left (\frac{ a_l -\sqrt{\Delta_l }}{a_l^2+b_l^2}, \frac{a_l  + \sqrt{\Delta_l}}{a_l^2+b_l^2}
\right )$, the inequality (\ref{minus}) is satisfied.
Thus,  there exists real $\mu$ such that the matrix $(I_N\otimes S_0-\mu(H\otimes S_0))$ is Schur  if and only if conditions (\ref{con1}) and (\ref{mu1})  are satisfied. And,
if conditions (\ref{con1}) and (\ref{mu1})  are satisfied, then, for all $\mu$ satisfying  (\ref{conmu}),
the matrix $(I_N\otimes S_0-\mu(H\otimes S_0))$ is Schur.
\end{Proof}
\end{Lemma}

By imposing some additional conditions on  the digraph $\mathcal{\bar{G}}$ or the matrix $S_0$, we can obtain two special cases of Lemma 3.1 as follows.

\begin{Corollary}

\textit{Case (i):}  Under Assumption \ref{Ass2.5}  and the additional assumption that the digraph $\mathcal{\bar{G}}$  is undirected,
the matrix $(I_N\otimes S_0-\mu(H\otimes S_0))$ is Schur if and only if
  \EQ \label{mu1x}
\max_{l = 1, \cdots, N} \bigg \{ \frac{|\lambda_q| -1}{a_l |\lambda_q|}  \bigg \} <   \min_{l = 1, \cdots, N}  \bigg \{  \frac{|\lambda_q| + 1}{a_l |\lambda_q|} \bigg \}.
\EN
And, if the condition (\ref{mu1x}) is satisfied, the matrix $(I_N\otimes S_0-\mu(H\otimes S_0))$ is Schur for all $\mu$ satisfying
\EQ \label{mu1xx}
\max_{l = 1, \cdots, N} \bigg \{ \frac{|\lambda_q| -1}{a_l |\lambda_q|}  \bigg \} <  \mu < \min_{l = 1, \cdots, N}  \bigg \{  \frac{|\lambda_q| + 1}{a_l |\lambda_q|} \bigg \}.
\EN

\textit{Case (ii):} Under Assumptions \ref{Ass2.5} and \ref{Ass2.7},
the matrix $(I_N\otimes S_0-\mu(H\otimes S_0))$ is Schur for all $\mu$ satisfying  (\ref{conmu}).


\begin{Proof}

 \textit{Case (i):}  By Remark \ref{rem23.1}, in this case, the matrix $H$ is symmetric and positive definite and hence $b_l = 0$ for $l=1,\cdots,N$. Thus,
the condition (\ref{con1}) is satisfied automatically, and the condition (\ref{mu1}) simplifies to (\ref{mu1x}) since, for $ l=1,\cdots,N$, $\Delta_{l} = \frac{ a_l^2}{|\lambda_q|^2 }$.


\textit{Case (ii):}
If  $|\lambda_q | \leq 1$, then the condition (\ref{con1}) is satisfied obviously, and the condition \eqref{mu1} is also satisfied since $\Delta_l \geq a_l$.

\end{Proof}

\end{Corollary}

\begin{Remark}
Under Assumptions \ref{Ass2.5}, and  \ref{Ass2.7},
the observer (\ref{ob2}) is always asymptotically stable for some real $\mu$. Nevertheless, this is not the case for the
observer (\ref{ob}). In fact, denote the eigenvalues of $S_0$ by $\{\alpha_k\pm j\beta_k\}$, where $\beta_k=0$ when $1\leq k \leq q_1$ with $0 \leq q_1 \leq q$, and $\beta_k\neq 0$ when $q_1<k\leq q_2$ where $q_1+2(q_2-q_1)=q$.
Then, by Theorem 1 in \cite{Su1}, the eigenvalues of $(I_N\otimes S_0-\mu(H\otimes I_q))$ are
\EQ \label{eig2}
\{\alpha_k\pm j\beta_k-\mu (a_j\pm j b_j):k=1,\cdots,q, j=1,\cdots,N\}
\EN
Thus the matrix  $(I_N\otimes S_0-\mu(H\otimes I_q))$ is Schur if and only if, for $k=1,\cdots,q$ and $j=1,\cdots,N$,
\EQQ
|\alpha_k\pm j\beta_k-\mu (a_j\pm j b_j)|<1
\ENN
and if and only if, for $k=1,\cdots,q$ and $j=1,\cdots,N$,
\EQ \label{eqleq12}
(\alpha_k-\mu a_j)^2+(\beta_k\pm\mu b_j)^2<1.
\EN
For each pair of $k$ and $j$, (\ref{eqleq12}) is equivalent to the following two inequalities:
\EQ \label{minus2}
(a_j^2+b_j^2)\mu^2-2\mu(a_j\alpha_k-b_j\beta_k)+\alpha_k^2+\beta_k^2-1<0
\EN
and
\EQ \label{plus}
 (a_j^2+b_j^2)\mu^2-2\mu(a_j\alpha_k+b_j\beta_k)+\alpha_k^2+\beta_k^2-1<0
\EN
Assume that the digraph $\mathcal{\bar{G}}$  is undirected and all the eigenvalues of $S_0$ have modulus $1$. Then it can be easily concluded that
the matrix $(I_N\otimes S_0-\mu(H\otimes I_q))$ is Schur for some real $\mu$  only if, for $k=1,\cdots,q$, $\alpha_k$ have the same sign.
Thus, as shown at the beginning of this Section,
if $S_0 = \mbox{diag} [1, -1]$, and the digraph $\mathcal{\bar{G}}$  is undirected, the  observer (\ref{ob}) cannot be asymptotically stable for any real $\mu$.
\end{Remark}

\section{Main Results}\label{cmr}

In this section, we will present our main result. For this purpose, we first present a simple result as follows.
\begin{Lemma}\label{Lem22}
Consider the system
\begin{equation}\label{sw1}
\begin{split}
&\zeta(t+1) = \sum_{i=0}^{h} F_i \zeta(t-\tau_i) + \sum_{i=0}^h G_i \xi(t-\tau_i)\\
&\xi(t+1) = \sum_{i=0}^h H_i \xi(t-\tau_i)\\
&\zeta(\theta)= \zeta_0(\theta),  \  \xi(\theta)=\xi_0(\theta), \ \theta \in I[-\tau,0]
\end{split}
\end{equation}
where $F_i \in \R^{n\times n}$, $G_i \in \R^{n \times m}$, $H_i \in \R^{m \times m},\ i=0,1,\dots,h$,
$0=\tau_0<\tau_1<\tau_2<\dots<\tau_h=\tau$ are non-negative integers.  Assume that $\zeta(t+1) = \sum_{i=0}^{h} F_i \zeta(t-\tau_i)$ and $\xi(t+1) = \sum_{i=0}^h H_i \xi(t-\tau_i)$ are exponentially stable. Then, the system \eqref{sw1} is exponentially stable.
\end{Lemma}
\begin{Proof}
Let $\eta(t)= \mbox{col} \left( \zeta(t), \xi(t)\right)$. Then the system \eqref{sw1} becomes
\begin{equation}\label{lem3.3sys}
\eta (t+1)=\sum_{i=0}^h \left( \begin{array}{cc}
F_i & G_i\\
0 & H_i\\
\end{array} \right) \eta(t-\tau_i)
\end{equation}

Let $M_i=\left( \begin{array}{cc}
F_i & G_i\\
0 & H_i\\
\end{array}\right), i=0,1,\dots, h$. Then the characteristic function of \eqref{lem3.3sys} is
\begin{equation}\label{root}
\begin{aligned}
\Delta(z) =& det{\left( z I_{n+m}-\sum_{i=0}^h M_i z^{-\tau_i} \right)}\\
=& det{\bigg(
\begin{array}{cc}
z I_n-\sum_{i=0}^h F_iz^{-\tau_i} & -\sum_{i=0}^h G_iz^{-\tau_i}\\
0& z I_m-\sum_{i=0}^h H_iz^{-\tau_i}\\
\end{array}
\bigg)}\\
=& det{\left(z I_n-\sum_{i=0}^h F_iz^{-\tau_i}\right)}det{\left(z I_m-\sum_{i=0}^h H_iz^{-\tau_i}\right)}
\end{aligned}
\end{equation}
Since $\zeta(t+1) = \sum_{i=0}^{h} F_i \zeta(t-\tau_i)$ and $\xi(t+1) = \sum_{i=0}^h H_i \xi(t-\tau_i)$ are exponentially stable, by Remark \ref{rem23.1}, all the roots of $\Delta(z)$ have modulus smaller than $1$. Thus the system \eqref{sw1} is exponentially stable.
\end{Proof}

\begin{Theorem}\label{Them2}
Under Assumptions \ref{Ass2.2}-\ref{Ass2.5}, if there exits a $\mu$ such that the distributed observer (\ref{obt}) is exponentially stable, then the cooperative output
regulation problem of the multi-agent system composed of \eqref{sys1}, \eqref{exo1}, \eqref{distur} is solvable by the distributed dynamic
output feedback control law of the form \eqref{ctr22}.
\end{Theorem}

\begin{Proof}
For $i=1,\cdots,N$, let $(X_i,\ U_i)$ satisfy the discrete regulator equations \eqref{re23}, $K_{1i} \in \R^{m_i\times n_i}$ satisfy Assumption \ref{Ass2.2}, and $K_{2i}=U_i-K_{1i}X_i$. Partition $K_{2i}$ as $K_{2i}=[K_{2xi},\ K_{2wi}]$ where $K_{2xi} \in \R^{m_i\times q}$ and $K_{2wi} \in \R^{m_i\times s_i}$. Define
\begin{equation}\label{control}
\begin{split}
u_i(t)=&\left[\begin{array}{cc}
            K_{1i}&K_{2wi}\\
            \end{array} \right]\left[\begin{array}{cc}
\hat{x}_i(t)\\
\hat{w}_i(t)
\end{array}\right]+K_{2xi}\eta_i(t)\\
\left[\begin{array}{cc}
\hat{x}_i(t+1)\\
\hat{w}_i(t+1)
\end{array}\right]=&\left[\begin{array}{cc}
        A_i&E_{wi}\\
        0&Q_i\\
        \end{array}\right]\left[\begin{array}{cc}
\hat{x}_i(t)\\
\hat{w}_i(t)
\end{array}\right]+\sum_{l=0}^h\left[\begin{array}{cc}
                        B_{il}\\
                        0\\
                        \end{array}\right]u_i(t-r_l)
                        +\left[\begin{array}{cc}
                        E_{xi}\\
                        0\\
                        \end{array}\right]\eta_i(t)\\&+L_i\bigg[y_{mi}(t)-\left(\begin{array}{cc}
                        \bar{C}_{i}&\bar{F}_{wi}\\
                        \end{array}\right)
                        \left[\begin{array}{cc}
\hat{x}_i(t)\\
\hat{w}_i(t)
\end{array}\right]
        -\sum_{l=0}^h\bar{D}_{il}u_i(t-r_l)-\bar{F}_{xi}\eta_i(t)\bigg]\\
\end{split}
\end{equation}
where  $L_i$ is such that the matrix $\bigg(\left[\begin{array}{cc}
                                        A_i&E_{wi}\\
                                        0&Q_i
                                        \end{array}\right]-L_i\left[\begin{array}{cc}\bar{C}_i&\bar{F}_{wi}
                                        \end{array}\right]\bigg)$ is Schur which exists under Assumption \ref{Ass2.3}.

Our control law is composed of  (\ref{control}) and the distributed observer (\ref{obt}) which can be verified to be in the form \eqref{ctr22}. We will show that under this control law, the closed-loop system satisfies the
two properties described in Definition 2.1. For this purpose, for $i=1,\cdots,N$, let $K_{\xi i} = [K_{1i},\ K_{2wi}]$, $\xi_i = \mbox{col}(\hat{x}_i, \hat{w}_i)$,
$v_i=\mbox{col}(x_0, w_i)$, $\bar{x}_i=x_i-X_iv_i$, $\bar{u}_i=u_i-U_iv_i$, $\xi_{ei}=\xi_i-\mbox{col}(x_i, \ w_i)$, and $\tilde{\eta}_i=\eta_i-x_0$. Then it can be verified that


\begin{equation}\label{ui}
\begin{split}
\bar{u}_i=&\left[\begin{array}{cc}
K_{1i}&K_{2wi}
\end{array}\right]\xi_i+K_{2xi}\eta_i-U_iv_i\\
=&\left[\begin{array}{cc}
K_{1i}&K_{2wi}
\end{array}\right]\left[\begin{array}{cc}
                        x_i\\
                        w_i
                        \end{array}\right]+K_{\xi i}\xi_{ei}+K_{2xi}x_0+K_{2xi}\tilde{\eta}_i-U_iv_i\\
                        =&K_{1i}x_i+K_{2wi}w_i+K_{2xi}x_0+K_{\xi i}\xi_{ei}+K_{2xi}\tilde{\eta}_i-(K_{2i}+K_{1i}X_i)v_i\\
                        =&K_{1i}(x_i-X_iv_i)+K_{2i}v_i+K_{\xi i}\xi_{ei}+K_{2i}\tilde{\eta}_i-K_{2i}v_i\\
                        =&K_{1i}\bar{x}_i+K_{\xi i}\xi_{ei}+K_{2xi}\tilde{\eta}_i\\
\end{split}
\end{equation}

\begin{equation}\label{xi}
\begin{split}
\bar{x}_i(t+r+1)=&x_i(t+r+1)-X_iv_i(t+r+1)\\
=&A_ix_i(t+r)+\sum_{l=0}^hB_{il}u_i(t+r-r_l)+E_iv_i(t+r)-X_i\bar{S}_i^{r+1}v_i(t)\\
=&A_i\big(\bar{x}_i(t+r)+X_iv_i(t+r)\big)+\sum_{l=0}^hB_{il}\big(\bar{u}_i(t+r-r_l)+U_iv_i(t+r-r_l)\big)\\&+E_iv_i(t+r)-X_i\bar{S}_i^{r+1}v_i(t)\\
=&A_i\bar{x}_i(t+r)+\sum_{l=0}^hB_{il}\bar{u}_i(t+r-r_l)+\big(A_iX_i\bar{S}_i^r+\sum_{l=0}^hB_{il}U_i\bar{S}_i^{r-r_l}+E_i\bar{S}^r_i\\&-X_i\bar{S}_i^{r+1}\big)v_i(t)\\
=&A_i\bar{x}_i(t+r)+\sum_{l=0}^hB_{il}\bar{u}_i(t+r-r_l)\\
\end{split}
\end{equation}

 \begin{equation} \label{zei}
\begin{split}
\xi_{ei}(t+1)=&\xi_i(t+1)-\left[\begin{array}{cc}
                    x_i(t+1)\\
                    w_i(t+1)\\
                    \end{array}\right]\\
                    =&\bigg(\left[\begin{array}{cc}
                                        A_i&E_{wi}\\
                                        0&Q_i
                                        \end{array}\right]-L_i\left[\begin{array}{cc}\bar{C}_i&\bar{F}_{wi}
                                        \end{array}\right]\bigg)\xi_{ei}+\bigg(\left[\begin{array}{cc}
                                        E_{xi}\\
                                        0\\
                                        \end{array}\right]-L_i\bar{F}_{xi}\bigg)\tilde{\eta}_i\\
\end{split}
\end{equation}

\begin{equation} \label{eqei}
\begin{split}
e_i(t+r)=&C_ix_i(t+r)+\sum_{l=0}^hD_{il}u_i(t+r-r_l)+F_{xi}x_0(t+r)+F_{wi}w_i(t+r)\\
=&C_i\big(\bar{x}_i(t+r)+X_iv_i(t+r)\big)+\sum_{l=0}^hD_{il}\bar{u}_i(t+r-r_l)+\sum_{l=0}^hD_{il}U_iv_i(t+r-r_l)\\&+F_iv_i(t+r)\\
=&C_i\bar{x}_i(t+r)+\sum_{l=0}^hD_{il}\bar{u}_i(t+r-r_l)+\big(C_iX_i\bar{S}_i^{r}+\sum_{l=0}^hD_{il}U_i\bar{S}_i^{r-r_l}+F_i\bar{S}_i^{r}\big)v_i(t)\\
=&C_i\bar{x}_i(t+r)+\sum_{l=0}^hD_{il}\bar{u}_i(t+r-r_l)\\
\end{split}
\end{equation}
Substituting \eqref{ui} into \eqref{xi} gives
\begin{equation}
\begin{split}
\bar{x}_i(t+r+1)=&A_i\bar{x}_i(t+r)+\sum_{l=0}^hB_{il}K_{1i}\bar{x}_i(t+r-r_l)+\sum_{l=0}^hB_{il}K_{\xi i}\xi_{ei}(t+r-r_l)\\&+\sum_{l=0}^hB_{il}K_{2xi}\tilde{\eta}_i(t+r-r_l)\\
\end{split}
\end{equation}

Under Assumption \ref{Ass2.2}, system $\bar{x}_i(t+1)=A_i\bar{x}_i+\sum_{l=0}^hB_{il}K_{1i}\bar{x}_i(t-r_l)$ is exponentially stable, and
under Assumption \ref{Ass2.3}, we can find $L_i$ such that system $\xi_{ei}(t+1)=
                    \bigg(\left[\begin{array}{cc}
                                        A_i&E_{wi}\\
                                        0&Q_i
                                        \end{array}\right]-L_i\left[\begin{array}{cc}\bar{C}_i&\bar{F}_{wi}
                                        \end{array}\right]\bigg)\xi_{ei}(t)$
is exponentially stable. Moreover, by the assumption of the Theorem, the distributed observer (\ref{obt}) is exponentially stable. Thus,
by Lemma \ref{Lem22}, the closed-loop system composed of (\ref{xi}), (\ref{zei}) and the distributed observer (\ref{obt}) is exponentially stable.
Hence, we have, for $i=1,\cdots, N$, $\lim_{t\rightarrow \infty} \bar{x}_i(t)=0$,  $\lim_{t\rightarrow \infty}\xi_{ei}(t)=0$, and $\lim_{t\rightarrow \infty}\tilde{\eta}_i=0$.
Therefore, from \eqref{ui}, we have, for $i=1,\cdots, N$, $\lim_{t\rightarrow \infty}\bar{u}_i(t)=0$. Finally, from (\ref{eqei}), we have,
for $i=1,\cdots, N$, $\lim_{t\rightarrow \infty}e_i(t)=0$. Thus the proof is complete.
\end{Proof}

\begin{Remark}
Under the additional Assumption \ref{Ass2.7},
  by \textit{Case (ii)} of Corollary 3.1,
there always exists some real $\mu$ satisfying (\ref{conmu}) such that
 the distributed observer (\ref{obt}) is asymptotically stable.
In this case,  the  cooperative output
regulation problem of the multi-agent system composed of \eqref{sys1}, \eqref{exo1} and \eqref{distur} is always solvable by the distributed dynamic
output feedback control law of the form (\ref{ctr22}).
\end{Remark}

\begin{Remark}
The delay-free system can be viewed as a special case of the system \eqref{sys1} with $h=0$. In this case, Assumption \ref{Ass2.2} reduces to the following:
\end{Remark}

\begin{Assumption}\label{Ass2.8}
For $i = 1, \cdots, N$, $(A_i, B_{i0})$ is stabilizable.
\end{Assumption}

For convenience, we state the solvability of this special case as follows:

\begin{Theorem}\label{Them3}
Under Assumptions \ref{Ass2.3}-\ref{Ass2.5}, and \ref{Ass2.8}, if there exits a $\mu$ such that the distributed observer (\ref{obt}) is exponentially stable, then the cooperative output
regulation problem of the multi-agent system composed of \eqref{sys1} with $h=0$, \eqref{exo1}, and \eqref{distur} is solvable by the distributed dynamic
output feedback control law of the  following form:

\begin{equation}\label{controldelayfree}
\begin{split}
u_i(t)= &\left[\begin{array}{cc}
            K_{1i}&K_{2wi}\\
            \end{array} \right]\left[\begin{array}{cc}
\hat{x}_i(t)\\
\hat{w}_i(t)
\end{array}\right]+K_{2xi}\eta_i(t)\\
            \left[\begin{array}{cc}
\hat{x}_i(t+1)\\
\hat{w}_i(t+1)
\end{array}\right] = &\left[\begin{array}{cc}
        A_i&E_{wi}\\
        0&Q_i\\
        \end{array}\right]\left[\begin{array}{cc}
\hat{x}_i(t)\\
\hat{w}_i(t)
\end{array}\right]+ \left[\begin{array}{cc}
                        B_{i}\\
                        0\\
                        \end{array}\right] u_i(t)
                        +\left[\begin{array}{cc}
                        E_{xi}\\
                        0\\
                        \end{array}\right]\eta_i(t)+\\
                        &L_i\bigg (y_{mi}(t)- [
                        \bar{C}_{i} ~ \bar{F}_{wi} ]\left[\begin{array}{cc}
\hat{x}_i(t)\\
\hat{w}_i(t)
\end{array}\right]
        - \bar{D}_{i}u_i(t)-\bar{F}_{xi}\eta_i(t)\bigg )\\
        \eta_i(t+1)=&S_0\eta_i+\mu S_0\bigg(\sum_{j\in\bar{\mathcal{N}}_i}a_{ij}\big(\eta_j(t)-\eta_i(t)\big)\bigg)\\
        i = &1. \cdots, N
\end{split}
\end{equation}

\end{Theorem}

Needless to say that, under the additional Assumption \ref{Ass2.7}, a control law of the form (\ref{controldelayfree}) that solves the delay-free case always exists.

\section{An Example}\label{ste}
Consider the discrete time-delay multi-agent systems of the form \eqref{sys1}
with $N = 4$, $h=1, \ r_1=1$,
$A_i=\left[\begin{array}{cc}
                    1&1\\
                    0&1\\
                    \end{array}\right]$, $B_{i0}=\left[\begin{array}{cc}
                                            -0.5\\
                                            0
                                            \end{array}\right]$, $B_{i1}=\left[\begin{array}{cc}
                            1\\
                            1
                            \end{array}\right]$, $E_{xi}=\left[\begin{array}{cc}
                            -2\cos{1}+1.5&-1\\
                            -\cos{1}-\sin{1}&\cos{1}-1+\sin{1}
                            \end{array}\right]$, $\bar{C}_i=C_i=\left[\begin{array}{cc}
                            1&0
                            \end{array}\right]$, $\bar{D}_i=D_i=0$ ,  $\bar{F}_{xi}=F_{xi}=\left[\begin{array}{cc}
     1&0
     \end{array}\right]$, $\bar{F}_{wi}=F_{wi}=\left[\begin{array}{cc}
     -1&0
     \end{array}\right]$, $i=1,\cdots,4$, and
\begin{equation*}
\begin{split}
E_{w1}=&\left[\begin{array}{cc}
             -1.5+2\cos{2}+0.5\sin{2}&-1.25+0.5\cos{2}\\
     -0.5\sin{2}+\cos{2}&1.5\cos{2}-1-\sin{2}
     \end{array}\right],\\
E_{w2}=&\left[\begin{array}{cc}
             -1.5+2\cos{3}+\sin{3}&-1.5+\cos{3}\\
     \cos{3}&2\cos{3}-1-\sin{3}
     \end{array}\right],\\
E_{w3}=&\left[\begin{array}{cc}
             -1.5+2\cos{4}+1.5\sin{4}&-1.75+1.5\cos{4}\\
     0.5\sin{4}+\cos{4}&2.5\cos{4}-1-\sin{4}
     \end{array}\right],\\
E_{w4}=&\left[\begin{array}{cc}
             -1.5+2\cos{5}+2\sin{5}&-2+2\cos{5}\\
     \sin{5}+\cos{5}&3\cos{5}-1-\sin{5}
     \end{array}\right].\\
     \end{split}
     \end{equation*}

Let the  leader system be
\begin{equation*}
\begin{split}
x_0(t+1)&=S_0 x_0(t)\\
\end{split}
\end{equation*}
where $S_0=\left[ \begin{array}{cc}
\cos{1}&\sin{1}\\
-\sin{1}&\cos{1}
\end{array}\right]$.

The disturbances to four followers are generated by \eqref{distur}
with $Q_1=\left[ \begin{array}{cc}
\cos{2}&\sin{2}\\
-\sin{2}&\cos{2}
\end{array}\right]$, $Q_2=\left[ \begin{array}{cc}
\cos{3}&\sin{3}\\
-\sin{3}&\cos{3}
\end{array}\right]$, $Q_3=\left[ \begin{array}{cc}
\cos{4}&\sin{4}\\
-\sin{4}&\cos{4}
\end{array}\right]$, $Q_4=\left[ \begin{array}{cc}
\cos{5}&\sin{5}\\
-\sin{5}&\cos{5}
\end{array}\right].$

Letting $K_{1i}=\left(
                      \begin{array}{cc}
                        -0.0750 & -0.4650 \\
                      \end{array}
                    \right)$ for $i=1,\cdots,4,$ gives the roots of the polynomial $\mbox{det} \left (z_i I_2 - A_i - B_{i0} K_{1i} - B_{i1} K_{1i} z_i^{-1} \right)$, as $$\{0.6435 \pm 0.4436j,  0.7530\}. $$
                    Thus, Assumption \ref{Ass2.2} is satisfied. Also, it can be verified that Assumption \ref{Ass2.3} is satisfied.
Solving the regulator equations gives $X_i=\left[\begin{array}{cccc}
                    -1&0&1&0\\
                    0&1&0&1
                    \end{array}\right]$ and $U_i=\left[
                    \begin{array}{cccc}
                     1&0&-1& -0.5i \\
                      \end{array}
                      \right]
$. Thus, Assumption \ref{Ass2.4} is also satisfied.

\begin{figure}[H]
\centering
\includegraphics[scale=0.25]{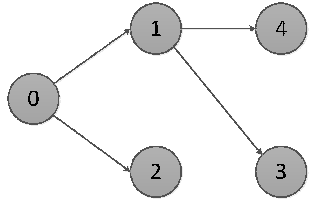}
\caption{Network topology} \label{network2}
\end{figure}

The network topology of the five agents is described in Fig. \ref{network2} which satisfies Assumption  \ref{Ass2.5}.

Moreover, the leader's signal is a sinusoidal function. By \textit{Case (ii)} of Corollary 3.1, the distributed observer for the leader system always exists.
In fact, it can be seen that
\begin{equation*}
\begin{split}
H=\left[ \begin{array}{cccc}
            1&0&0&0\\
            0&1&0&0\\
            -1&0&1&0\\
            -1&0&0&1\\
            \end{array}\right]
\end{split}
\end{equation*}
\normalsize
whose eigenvalues  are $\{1,1,1,1\}$. Thus, from (\ref{conmu}), for all $\mu$ satisfying $0<\mu<2$, the matrix $(I_N\otimes S_0-\mu(H\otimes S_0))$ is Schur.
Let us take $\mu = 0.5$.

Thus, by Theorem \ref{Them2}, the cooperative output regulation
problem for this example is solvable. For $i=1,\cdots,4$, with $K_{1i}=\left[\begin{array}{cc}
-0.075&-0.465\\
\end{array}\right]$, we can obtain $K_{2i}=U_i-K_{1i}X_i= \left[
                  \begin{array}{cccc}
   0.925&0.465&-0.925&-0.5i+0.465\\
                  \end{array}
                \right]
$. Thus,
$K_{2wi}=\left[\begin{array}{cc}
-0.925&-0.5i+0.465\\
\end{array}\right]$, $K_{2xi}=\left[\begin{array}{cc}
0.925&0.465\\
\end{array}\right]$.\\
 Finally,  let
\begin{equation*}
\begin{split}
L_{1}&=\left[\begin{array}{cccc}
0.0293&-1.5213&0.2372&-0.9982\\
\end{array}\right]^T\\
 L_2&=\left[\begin{array}{cccc}
-10.5975&3.6208&-9.2420&4.4453\\
\end{array}\right]^T\\
L_3&=\left[\begin{array}{cccc}
-1.6276&-1.1082&-0.9447&-0.4423\\
\end{array}\right]^T\\
L_4&=\left[\begin{array}{cccc}
1.3550&0.1854&0.1633&0.0026\\
\end{array}\right]^T
\end{split}
\end{equation*}
which make $\bigg(\left[\begin{array}{cc}
                                        A_i&E_{wi}\\
                                        0&Q_i
                                        \end{array}\right]-L_i\left[\begin{array}{cc}\bar{C}_i&\bar{F}_{wi}
                                        \end{array}\right]\bigg)$, $i = 1,2,3,4$,  Schur  matrices.

\begin{figure}[H]
\centering
\includegraphics[scale=0.6]{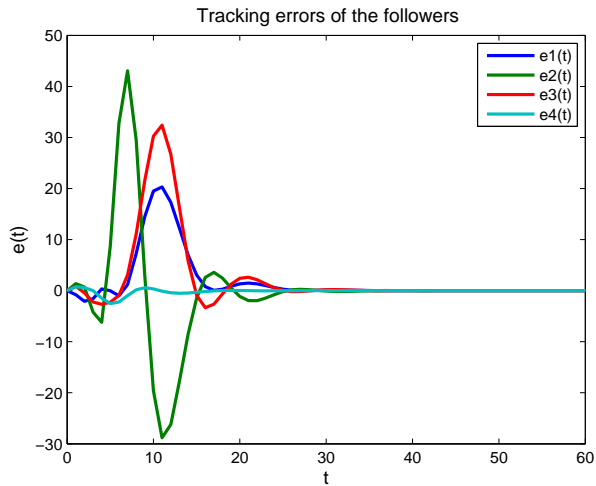}
\caption{Tracking errors of the followers} \label{error3}
\end{figure}

\begin{figure}[H]
\centering
\includegraphics[scale=0.6]{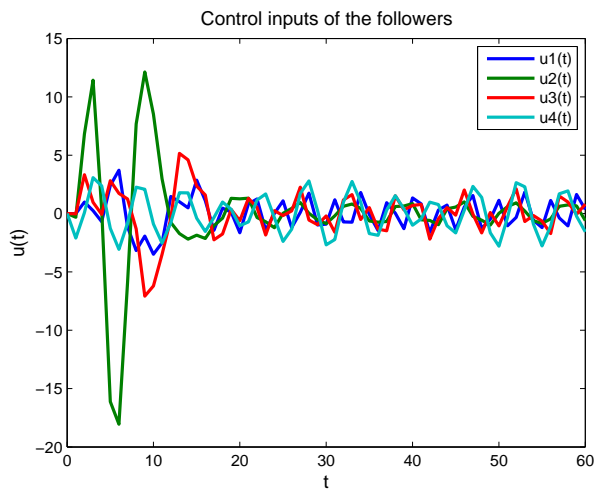}
\caption{Control inputs of the followers} \label{controlinput}
\end{figure}

Simulation is conducted for random initial conditions.  Fig. \ref{error3} shows the tracking performance  of the  error output $e_i$ of the fours followers and it can be seen that the distributed control law solves the cooperative output regulation problem successfully.
Fig. \ref{controlinput} further shows the controls of the four followers.


\section{Conclusion}\label{cc}
In this paper, we have studied the cooperative output regulation problem of discrete-time
linear time-delay multi-agent systems  under static network topology.
We have first thoroughly
studied the existence condition of the discrete distributed observer, and then
presented the solvability of the problem by distributed dynamic output feedback control law.
A natural consideration for our future work is the  study of the same problem for the same systems subject to the switching network topology.

\begin{appendix}
\subsection{Graph}
 A digraph $\mathcal{G}=(\mathcal{V},\mathcal{E})$ consists of a node set $\mathcal{V}=\{1,\cdots,N\}$ and an edge set $\mathcal{E}\subseteq \mathcal{V}\times \mathcal{V}$. An edge of $\mathcal{E}$ from node $i$ to node $j$ is denoted by $(i,j)$, where the nodes $i$ and $j$ are called the parent node and the child node of each other, and the node $i$ is also called a neighbor of the node $j$. Let $\mathcal{N}_i=\{j,(j,i) \in \mathcal{E}\}$ denote the subset of $\mathcal{V}$ which consists of all the neighbors of the node $i$. Edge $(i,j)$ is called undirected if $(i,j)\in \mathcal{E}$ implies that $(j,i)\in \mathcal{E}$. The graph is called undirected if every edge in $\mathcal{E}$ is undirected. If there exists a set of edges $\{(i_1,i_2),\cdots,(i_k,i_{k+1})\}$ in the digraph $\mathcal{G}$, then $i_{k+1}$ is said to be reachable from node $i_1$. A digraph $\mathcal{G}_s=(\mathcal{V}_s,\mathcal{E}_s)$, where $\mathcal{V}_s\subseteq \mathcal{V}$ and $\mathcal{E}_s \subseteq \mathcal{E}\cap (\mathcal{V}_s \times \mathcal{V}_s)$, is a subgraph of the digraph $\mathcal{G}=(\mathcal{V},\mathcal{E})$.
A weighted adjacency matrix of $\mathcal{G}$ is a square matrix denoted by $\mathcal{A}=\left[a_{ij}\right] \in R^{N\times N}$ such that, for $i,j = 1, \cdots, N$,  $a_{ii}=0$, $a_{ij}>0$ $\Leftrightarrow$ $(j,i)\in \mathcal{E}$, and $a_{ij} = a_{ji}$ if $(i,j)$  is undirected.
The Laplacian matrix of a digraph $\mathcal{G}$ is denoted by $\mathcal{L}=[l_{ij}]\in R^{N\times N}$, where $l_{ii}=\sum_{j=1}^Na_{ij}$, and $l_{ij}=-a_{ij}$ if $i\neq j$. More detailed exposition on graph theory can be found in \cite{Graph}.
\end{appendix}

\end{document}